\documentclass[preprint,11pt,a4paper,nopreprintline]{elsarticle}
\usepackage[T1]{fontenc}
\usepackage{lmodern}
\usepackage{amsmath,amssymb,amsthm}
\usepackage{graphicx}
\graphicspath{{figures/}{./}}
\usepackage{microtype}
\usepackage{etoolbox}
\usepackage[hidelinks]{hyperref}
\biboptions{numbers,sort&compress}

\journal{Advances in Mathematics}

\makeatletter
\patchcmd{\pprintMaketitle}{\hrule\vskip12pt}{\hrule width0pt\vskip12pt}{}%
  {\PackageError{manuscript}{Unable to remove the upper abstract rule}{}}
\patchcmd{\pprintMaketitle}{\hrule\vskip12pt}{\hrule width0pt\vskip12pt}{}%
  {\PackageError{manuscript}{Unable to remove the lower abstract rule}{}}
\makeatother

\numberwithin{equation}{section}
\newtheorem{theorem}{Theorem}[section]
\newtheorem{lemma}[theorem]{Lemma}
\newtheorem{proposition}[theorem]{Proposition}
\newtheorem{corollary}[theorem]{Corollary}
\theoremstyle{remark}
\newtheorem{remark}[theorem]{Remark}
\newtheorem{example}[theorem]{Example}
\theoremstyle{plain}
\newtheorem*{conjecture}{$\boldsymbol{k}$-harmonic conjecture}
\newcommand{\R}{\mathbb R}
\newcommand{\Sph}{\mathbb S}
\newcommand{\dd}{\,\mathrm d}
\newcommand{\vol}{\operatorname{Vol}}
\newcommand{\II}{\mathrm{II}}
\newcommand{\Hm}{\mathbf H}
\DeclareMathOperator{\tr}{tr}
\DeclareMathOperator{\Scal}{Scal}
\DeclareMathOperator{\diver}{div}
\hypersetup{
  pdftitle={Counterexamples to the k-harmonic conjecture},
  pdfauthor={Yu Fu, Yanru Luo, Dan Yang},
  pdfsubject={Biharmonic and polyharmonic hypersurfaces in Euclidean space},
  pdfkeywords={k-harmonic conjecture, biharmonic hypersurface, triharmonic hypersurface, circular cone}
}
\begin{document}
\begin{frontmatter}

\title{Counterexamples to the \texorpdfstring{$k$}{k}-harmonic conjecture}

\author[dufe-data]{Yu Fu}
\ead{yufu@dufe.edu.cn}
\author[dufe-data]{Yanru Luo}
\ead{luoyanru0819@163.com}
\author[liaoning]{Dan Yang}
\ead{dlutyangdan@126.com}

\affiliation[dufe-data]{
  organization={School of Data Science and Artificial Intelligence,
    Dongbei University of Finance and Economics},
  city={Dalian}, postcode={116025}, country={P.~R.~China}}
\affiliation[liaoning]{
  organization={School of Mathematics, Liaoning University},
  city={Shenyang}, postcode={110036}, country={P.~R.~China}}

\begin{abstract}
We construct nonminimal $k$-harmonic hypersurfaces in $\R^{n+1}$ for every
$k\geq3$ and $n\geq3$, giving counterexamples to the $k$-harmonic
conjecture. Each example has exact order $k$ and is an incomplete circular
cone or a Euclidean product of such a cone. We also classify the
$k$-harmonic circular cones, determine the possible dimensions of
nonminimal cones with closed links, and give constructions in higher
codimension.
\end{abstract}


\end{frontmatter}

\section{Introduction}

Let $\psi:(M^n,g)\to\R^N$ be an isometric immersion. Denote by $\II$ and
$\Hm=n^{-1}\tr\II$ the second fundamental form and the mean curvature
vector, respectively. We use the Laplacian $\Delta=-\diver\nabla$.
The Beltrami formula gives
\[
\Delta\psi=-n\Hm.
\]
For an integer $k\geq2$, the immersion is called $k$-harmonic if
\begin{equation}\label{eq:polyharmonic}
\Delta^k\psi=0,
\qquad\text{or equivalently}\qquad
\Delta^{k-1}\Hm=0.
\end{equation}
The Laplacian in \eqref{eq:polyharmonic} acts componentwise in the fixed
Euclidean coordinates and is computed with respect to the induced metric.
In particular, every minimal immersion is trivially $k$-harmonic.

The case $k=2$ is the equation of biharmonic submanifolds. Jiang \cite{Jiang1986}
established the first and second variational formulas for biharmonic maps.
In Euclidean space, the study of biharmonic submanifolds also arose from
Chen's theory of finite type submanifolds. Chen \cite{Chen1991}
conjectured that {\em every biharmonic submanifold of Euclidean space is
minimal}. The case of surfaces in $\R^3$ was proved independently by
Chen and Jiang \cite{Chen1991,Jiang1987}. Hasanis-Vlachos
\cite{HasanisVlachos1995} proved the case of hypersurfaces in $\R^4$,
and Fu-Hong-Zhan
\cite{FHZ2021,FHZ2023} subsequently settled the hypersurface cases
in $\R^5$ and $\R^6$. In a recent preprint \cite{FHY2026},
Fu-Hong-Yang prove that every biharmonic hypersurface in a space
form of nonpositive sectional curvature is minimal, in arbitrary
dimension.
In arbitrary codimension, Akutagawa and Maeta
\cite{AkutagawaMaeta2013} proved that every properly immersed biharmonic
submanifold of Euclidean space is minimal. For proper biharmonic hypersurfaces in spheres, there we refer to
\cite{BMO2008,BMO2010}. A closely
related higher-dimensional theory for the biharmonic hypersurface flow was developed by Fu, Hong and Tian \cite{FuHongTianJEMS} very recently.

Concerning  higher-order harmonic submanifolds, Maeta \cite{Maeta2012}
proposed the following extension of Chen's conjecture:

\begin{conjecture}
For each integer $k\geq2$, every $k$-harmonic submanifold of Euclidean
space is minimal.
\end{conjecture}

This formulation imposes no completeness or global integrability
assumptions. The examples below are incomplete, so they do not address
variants that require completeness. We use compactly supported
variations with the induced metric held fixed, as recalled in
Section~\ref{sec:main}.

Higher-order harmonic hypersurfaces have been studied under both global
and local geometric assumptions. Maeta, Nakauchi and Urakawa
\cite{MNU2015}, and Nakauchi and Urakawa \cite{NU2018}, obtained
minimality results under completeness and integrability assumptions.
Montaldo, Oniciuc and Ratto \cite{MOR2022} proved that a
$k$-harmonic hypersurface, $k\geq3$,
in a nonpositively curved space form is minimal if the mean curvature
and the squared norm of the second fundamental form are constant.
Chen and Guan \cite{ChenGuan2023} obtained rigidity results for
constant mean curvature (CMC) triharmonic hypersurfaces in
5-dimensional space forms.
Fu and Yang \cite{FuYang2023} proved minimality for CMC
triharmonic hypersurfaces in hyperbolic space and CMC triharmonic hypersurfaces in $\R^6$.

Branding et al. \cite{BMOR2020} studied the Eells--Sampson higher-order
energies and clarified their relationship with the usual $k$-energies.
For a Euclidean target, the two energies coincide, and both
Euler--Lagrange equations reduce to \eqref{eq:polyharmonic}.
Proper triharmonic immersions into spheres were constructed in
\cite{Maeta2015Construction}. Montaldo and Ratto \cite{MR2018}
constructed proper $k$-harmonic small hyperspheres and generalized
Clifford tori in spheres. Montaldo et al. \cite{MOR2024} classified
proper triharmonic isoparametric surfaces and CMC $k$-harmonic Hopf
cylinders, $k\geq3$, in three-dimensional
Bianchi--Cartan--Vranceanu spaces.
Further spherical examples and higher-order hypersurface equations
were studied in \cite{MOR2022}.
Branding \cite{Branding2023} investigated the normal stability of CMC
triharmonic hypersurfaces in space forms. Branding \cite{Branding2025}
also proved that the small hypersphere
$\Sph^n(1/2)\hookrightarrow\Sph^{n+1}$ has normal index one for both
the $4$-energy and the Eells--Sampson $4$-energy.

Ratto \cite{Ratto2020} derived the triharmonic equations for
$\mathrm{SO}(p+1)\times\mathrm{SO}(q+1)$-invariant hypersurfaces
in Euclidean space, including rotational hypersurfaces, and discussed
the compatibility of the resulting ordinary differential equations
with the isometric immersion condition. For the homogeneous cones
considered below, the induced metric gives an explicit factorization
of the iterated Laplacian. The direct calculation for a straight
profile is compared with the printed rotational equations in
Remark~\ref{rem:ratto}.

For comparison, the Bernstein problem asks whether every entire minimal graph in
$\R^{n+1}$ is a hyperplane. Bernstein himself \cite{Bernstein1915} proved
the case $n=2$, and Simons \cite{Simons1968} established the result
for $n\leq7$. Bombieri et al. \cite{BDG1969} constructed nonaffine
entire minimal graphs for $n\geq8$. The Bernstein problem concerns
flatness of entire minimal graphs, whereas the $k$-harmonic conjecture
concerns minimality of $k$-harmonic submanifolds.

Before stating the main results, we say that an immersion has \emph{exact order
$k$} if $\Delta^k\psi=0$ and $\Delta^{k-1}\psi\not\equiv0$.
An immersion of exact order $k$ is $p$-harmonic for every $p\geq k$.
Thus a nonminimal triharmonic immersion already gives counterexamples
for every $k\geq3$. Our first result prescribes the exact order:
for each $k\geq3$, we construct a family whose first vanishing
iterate is $\Delta^k\psi$. The examples are graphs over
$(\R^3\setminus\{0\})\times\R^{n-3}$ with incomplete induced metrics.

\begin{theorem}\label{thm:main}
Let $k\geq3$ and $n\geq3$ be integers, and put
$c_k=\sqrt{(k-1)(2k-5)}$. Then the hypersurface
\begin{equation}\label{eq:main-graph}
M_{k,n}
=
\bigl\{(x,z,c_k|x|):
x\in\R^3\setminus\{0\},\ z\in\R^{n-3}\bigr\}
\subset\R^{n+1}
\end{equation}
is a smooth embedded nonminimal $k$-harmonic hypersurface.
The inclusion $\psi_{k,n}$ has exact order $k$, and
$\Delta^{k-1}\psi_{k,n}$ is nowhere zero. The induced metric is incomplete.
\end{theorem}

For $n=3$, the first two examples are $x_4=\sqrt2\,|x|$ and
$x_4=3|x|$, of exact orders $3$ and $4$, respectively.
The case $k=3$ will be verified directly in
Example~\ref{ex:direct-triharmonic-cone}.

Our second result classifies the standard circular cones by dimension,
aperture and exact order. For $n\geq2$ and $0<a<1$, consider
\begin{equation}\label{eq:circular-cone}
\psi_a^n:(0,\infty)\times\Sph^{n-1}\longrightarrow\R^{n+1},
\qquad
\psi_a^n(s,\omega)=\bigl(as\omega,\sqrt{1-a^2}\,s\bigr).
\end{equation}

\begin{theorem}\label{thm:circular}
Let $n\geq2$ and $k\geq2$. The immersion $\psi_a^n$ has exact order $k$
if and only if
\begin{equation}\label{eq:aperture}
n\text{ is odd},\qquad 3\leq n\leq2k-3,
\qquad
a^2=\frac{n-1}{(2k-3)(2k-n-1)}.
\end{equation}
In particular, no nonplanar circular cone of even dimension is
$k$-harmonic for any finite $k$.
\end{theorem}

For a fixed $k\geq3$, the nonminimal $k$-harmonic standard circular cones are precisely those of exact order $j$, where $3\leq j\leq k$. For each such $j$, the cone dimension $n$ is odd and satisfies $3\leq n\leq2j-3$, and the aperture is given by \eqref{eq:aperture} with $k$ replaced by $j$.

 Our third result extends the dimension
restriction to arbitrary closed spherical links. For a nonminimal
$k$-harmonic cone, the spectral decomposition also forces the average
position of the link to be nonzero.

\begin{theorem}\label{thm:closed-link}
Let $k\geq2$ and $n\geq2$ be integers, and let
$\psi_L:(L^{n-1},h)\to\Sph^{N-1}(1)$ be an isometric immersion of a
connected closed manifold. Consider the cone
\[
\psi:(0,\infty)\times L\longrightarrow\R^N,
\qquad
\psi(s,x)=s\psi_L(x).
\]
If $\psi$ is nonminimal and $k$-harmonic, then $n$ is odd,
$3\leq n\leq2k-3$, and
\[
\int_L\psi_L\dd v_h\neq0.
\]
Conversely, for each odd $n$ with $3\leq n\leq2k-3$, there is a
nonminimal cone with closed link and exact order $k$ in $\R^{n+1}$.
\end{theorem}
\begin{remark}
	For $n=1$, every connected $k$-harmonic submanifold of Euclidean
	space is an open part of a straight line (see \cite{Maeta2012}).
	Consequently, no one-dimensional submanifold of Euclidean
	space has exact order $k\geq2$.
	
	For $n=2$, Theorems~\ref{thm:circular} and
	\ref{thm:closed-link} exclude nonminimal cones in the respective
	classes. The existence of a surface in $\R^3$ of exact order
	$k\geq3$ is not settled by the present results.
\end{remark}

The average position of a closed link also satisfies the following sharp
estimate.

\begin{corollary}\label{cor:centroid}
Let $\psi_L:(L^{n-1},h)\to\Sph^{N-1}(1)$ be an isometric immersion of a
connected closed manifold. Suppose that the cone $\psi(s,x)=s\psi_L(x)$
is nonminimal and $k$-harmonic, where $n\geq2$ and $k\geq2$. Put
\[
c=\frac1{\vol(L)}\int_L\psi_L\dd v_h,
\qquad
\Lambda=(2k-3)(2k-n-1).
\]
Then
\begin{equation}\label{eq:centroid-bound}
0<|c|^2\leq1-\frac{n-1}{\Lambda}.
\end{equation}
Equality in the upper bound holds if and only if, up to an orthogonal
transformation of $\R^N$,
\begin{equation}\label{eq:centroid-equality}
\psi_L(x)=\bigl(a\widehat\psi(x),\sqrt{1-a^2}\bigr),
\qquad
a^2=\frac{n-1}{\Lambda},
\end{equation}
where $0<a<1$ and $\widehat\psi:(L,a^{-2}h)\to\Sph^{N-2}(1)$ is a
minimal isometric immersion. In the equality case, the cone $\psi$ has exact
order $k$.
\end{corollary}

The proofs rely on the fact that the cone Laplacian lowers the degree of a
homogeneous function by two. We choose the aperture so that the
spherical component first vanishes at the prescribed iteration, while
the last component vanishes after two iterations. This gives the
three-dimensional examples, and Euclidean products yield
Theorem~\ref{thm:main}. For a closed link, the same iteration formula
and an integral identity restrict the possible nonzero spectral
components. The largest admissible spectral value gives the sharp
centroid bound, and equality reduces the link to a minimal immersion
into a small sphere.

We remark that Chen’s conjecture for biharmonic submanifolds of
arbitrary codimension in Euclidean space remains open. Although the
hypersurface case has been settled, the existence of nonminimal biharmonic
submanifolds in higher codimension is still unknown. Every nonminimal example constructed
in this paper has exact order $k\geq3$, and hence satisfies
$\Delta^2\psi\not\equiv0$.

The paper is organized as follows. In Section~\ref{sec:main} we give the
cone Laplacian formula and the proof of Theorem~\ref{thm:main}.
In Section~\ref{sec:classification}, we prove
Theorems~\ref{thm:circular} and \ref{thm:closed-link}, together with
Corollary~\ref{cor:centroid}.
Section~\ref{sec:examples} relates the equality case to constructions
from minimal spherical immersions and gives explicit examples in $\R^5$.

\section{Proof of Theorem~\ref{thm:main}}\label{sec:main}

We first recall the variational formulation. Fix the domain metric and
let $\Omega\subset M$ be a relatively compact smooth domain. For a
Euclidean-valued map, the $k$-energies are given by
(see \cite{Maeta2012,BMOR2020})
\[
\begin{aligned}
E_k(\psi;\Omega)
&=\frac12\int_\Omega|\Delta^{k/2}\psi|^2\dd v_g,
&& k\geq2\text{ even},\\[4pt]
E_k(\psi;\Omega)
&=\frac12\int_\Omega|\mathrm d\Delta^{(k-1)/2}\psi|^2\dd v_g,
&& k\geq1\text{ odd}.
\end{aligned}
\]
For any smooth variation field $V$ supported in $\Omega$, repeated
integration by parts yields
\[
\delta E_k(\psi;\Omega)[V]
=\int_\Omega\langle\Delta^k\psi,V\rangle\dd v_g.
\]
Hence \eqref{eq:polyharmonic} is the Euler--Lagrange equation for all
compactly supported variations. For an isometric immersion, the metric
is first taken to be the induced metric and is then fixed throughout
the variation.

For a hypersurface in Euclidean space $\R^{n+1}$, we choose a local
unit normal $\nu$ and write
\[
\Hm=H\nu,\qquad A(X)=-D_X\nu,\qquad nH=\tr A,
\]
where $D$ is the Euclidean connection and $H$ is the mean curvature.

The following lemma gives the iteration formula for the cone Laplacian.

\begin{lemma}\label{lem:cone}
Let $(L^{n-1},h)$ be a Riemannian manifold and equip
$(0,\infty)\times L$ with $g=\mathrm ds^2+s^2h$. For every smooth
function $v$ on $L$ and every $\alpha\in\R$, we have
\begin{equation}\label{eq:homogeneous}
\Delta_g(s^\alpha v)
=s^{\alpha-2}
\bigl(\Delta_hv-\alpha(\alpha+n-2)v\bigr).
\end{equation}
For integers $j\geq0$ and $p\geq1$, put
\begin{equation}\label{eq:cone-polynomial}
q_{j,n}=(1-2j)(n-1-2j),
\qquad
P_{p,n}(\lambda)=\prod_{j=0}^{p-1}(\lambda-q_{j,n}).
\end{equation}
Then
\begin{equation}\label{eq:iteration}
\Delta_g^p(sv)=s^{1-2p}P_{p,n}(\Delta_h)v.
\end{equation}
The same formulas hold componentwise for Euclidean-valued maps.
\end{lemma}

\begin{proof}
Since $\mathrm dv_g=s^{n-1}\mathrm ds\,\mathrm dv_h$, the cone Laplacian becomes
\[
\Delta_g
=-\partial_s^2-\frac{n-1}{s}\partial_s+\frac1{s^2}\Delta_h.
\]
Applying $\Delta_g$ to $s^\alpha v$, we obtain
\eqref{eq:homogeneous}. To prove \eqref{eq:iteration}, observe that
after $j$ iterations the radial degree is
$1-2j$, which gives the next factor
$\Delta_h-(1-2j)(n-1-2j)$. Multiplication over $j=0,\ldots,p-1$ yields
\eqref{eq:iteration}.
\end{proof}

\begin{figure}[htbp]
\centering
\includegraphics[width=\textwidth]{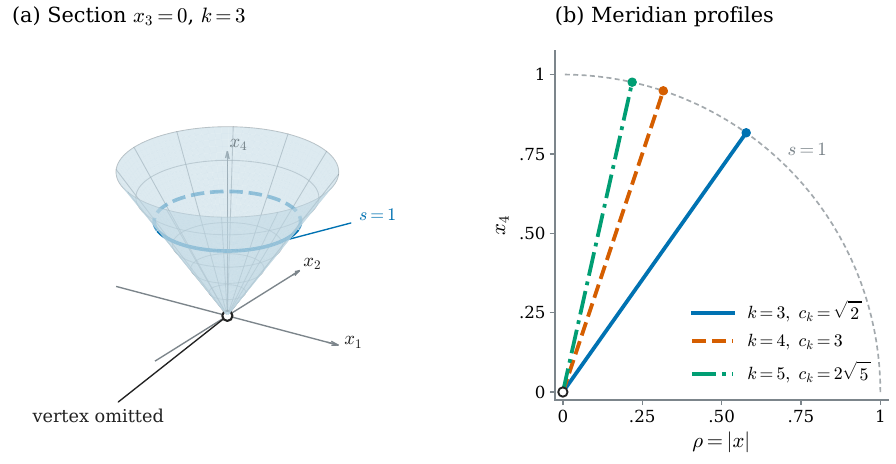}
\caption{Sections of $M_{k,3}\subset\R^4$.
(a) The slice $x_3=0$ for $k=3$, with the vertex removed and the
circle $s=1$ highlighted.
(b) Meridians $(\rho,x_4)=(a_ks,b_ks)$, $0<s\leq1$, for $k=3,4,5$.
Exact order $k$ refers to the full three-dimensional immersion.}
\label{fig:cone-sections}
\end{figure}

For later use, we record two identities which follow from
\eqref{eq:cone-polynomial}:
\begin{equation}\label{eq:root-difference}
\begin{aligned}
q_{j,n}-(n-1)&=2j(2j-n),\\
q_{j,n}-q_{\ell,n}
&=2(j-\ell)\bigl(2(j+\ell)-n\bigr).
\end{aligned}
\end{equation}
\begin{proof}[Proof of Theorem~\ref{thm:main}]
We first consider $n=3$. Set
\[
a_k=\frac1{\sqrt{(k-2)(2k-3)}},
\qquad
b_k=c_ka_k.
\]
Since
\[
(k-2)(2k-3)=1+(k-1)(2k-5),
\]
we have $a_k^2+b_k^2=1$. Define
\begin{equation}\label{eq:seed}
\psi_k(s,\omega)=(a_ks\omega,b_ks),
\qquad s>0,\quad\omega\in\Sph^2.
\end{equation}
For $V,W\in T_\omega\Sph^2$, we have
\[
\begin{aligned}
\mathrm d\psi_k(\partial_s)&=(a_k\omega,b_k),
&\mathrm d\psi_k(V)&=(a_ksV,0),\\
|\mathrm d\psi_k(\partial_s)|^2&=1,
&\langle\mathrm d\psi_k(\partial_s),\mathrm d\psi_k(V)\rangle&=0,\\
\langle\mathrm d\psi_k(V),\mathrm d\psi_k(W)\rangle
&=a_k^2s^2\langle V,W\rangle.
\end{aligned}
\]
It follows that the induced metric is given by
\begin{equation}\label{eq:seed-metric}
g_k=\mathrm ds^2+a_k^2s^2g_{\Sph^2}.
\end{equation}
Moreover, $x=a_ks\omega$ and $b_ks=c_k|x|$. Hence the image of
\eqref{eq:seed} is the smooth graph in \eqref{eq:main-graph} for $n=3$
(see Figure~\ref{fig:cone-sections}).

We next compute the iterated Laplacians. On the link
$(\Sph^2,a_k^2g_{\Sph^2})$, the spherical coordinate map satisfies
\[
\Delta_{a_k^2g_{\Sph^2}}\omega
=\lambda_k\omega,
\qquad
\lambda_k=\frac2{a_k^2}
=2(k-2)(2k-3)=q_{k-1,3}.
\]
The last component has link eigenvalue zero. Combining
Lemma~\ref{lem:cone} with \eqref{eq:seed-metric}, we obtain
\begin{equation}\label{eq:seed-iterates}
\Delta_{g_k}^p\psi_k
=s^{1-2p}
\bigl(a_kP_{p,3}(\lambda_k)\omega,\,
b_kP_{p,3}(0)\bigr).
\end{equation}
Since $q_{1,3}=0$, it follows that $P_{p,3}(0)=0$ for $p\geq2$.
Moreover, $\lambda_k=q_{k-1,3}$ yields
$P_{k,3}(\lambda_k)=0$, and hence $\Delta_{g_k}^k\psi_k=0$.

To verify that the order is exactly $k$, we use
\eqref{eq:root-difference}. For $0\leq j\leq k-2$, we obtain
\[
\lambda_k-q_{j,3}
=2(k-1-j)(2k+2j-5)>0.
\]
Consequently,
\[
C_k:=P_{k-1,3}(\lambda_k)
=\prod_{j=0}^{k-2}2(k-1-j)(2k+2j-5)>0.
\]
Using \eqref{eq:seed-iterates} again, we obtain
\begin{equation}\label{eq:last-iterate}
\Delta_{g_k}^{k-1}\psi_k
=\bigl(a_kC_ks^{3-2k}\omega,0\bigr)\neq0.
\end{equation}
Since $|\omega|=1$, the vector in \eqref{eq:last-iterate} is nowhere
zero, which shows that $\psi_k$ has exact order $k$.

We now compute the principal curvatures. Choose
$\nu=(-b_k\omega,a_k)$. Since $a_k^2+b_k^2=1$ and
$\langle\omega,V\rangle=0$, $\nu$ is a unit normal.
Differentiating $\nu$ yields
\[
\mathrm d\nu(\partial_s)=0,
\qquad
\mathrm d\nu(V)=(-b_kV,0)
=-\frac{b_k}{a_ks}\,\mathrm d\psi_k(V).
\]
Hence the Weingarten formula gives the three principal curvatures
\begin{align*}
0,\quad \frac{c_k}{s},\quad \frac{c_k}{s}.
\end{align*}
It follows that
\begin{equation}\label{eq:seed-curvature}
H=\frac{2c_k}{3s}>0,
\qquad
\Hm=\frac{2c_k}{3s}\nu,
\qquad
|A|^2=\frac{2c_k^2}{s^2}.
\end{equation}
Equation~\eqref{eq:seed-curvature} shows that $\psi_k$ is nonminimal
and has nonconstant mean curvature.
A radial curve has length element $\mathrm ds$ and reaches the missing
vertex in finite length, which implies that the induced metric is incomplete.

Finally, for $n>3$, after the ambient coordinate identification
$(u,t,z)\mapsto(u,z,t)$, the inclusion $\psi_{k,n}$ is parametrized by
$\psi_k\times\operatorname{id}_{\R^{n-3}}$.
By \eqref{eq:seed-metric}, the product metric is
$g_k+g_{\R^{n-3}}$. The product Laplacian is the sum of the factor
Laplacians. Since each coordinate function on $\R^{n-3}$ is harmonic,
under this identification we obtain
\[
\Delta^p\psi_{k,n}
=(\Delta_{g_k}^p\psi_k,0)
\qquad(p\geq1).
\]
Combining this identity with \eqref{eq:last-iterate}, we conclude that
$\psi_{k,n}$ has exact order $k$. The principal
curvatures are given by $c_k/s$ with multiplicity two and zero with multiplicity
$n-2$, respectively. Thus $H=2c_k/(ns)>0$. The metric remains incomplete, and the image is the graph in
\eqref{eq:main-graph}.
\end{proof}

\begin{example}\label{ex:direct-triharmonic-cone}
For $k=3$ and $n=3$, put $\rho=|x|=s/\sqrt3$.
In the coordinates $(\rho,\omega)$, the immersion $\psi_3$ in
\eqref{eq:seed} becomes
\[
\psi_3=(\rho\omega,\sqrt2\,\rho),
\qquad \rho>0,\quad \omega\in\Sph^2.
\]
The induced metric and its Laplacian are
\[
g=3\,\mathrm d\rho^2+\rho^2g_{\Sph^2},
\qquad
\Delta_g
=-\frac13\left(\partial_\rho^2+\frac2\rho\partial_\rho\right)
+\frac1{\rho^2}\Delta_{\Sph^2}.
\]
Since $\Delta_{\Sph^2}\omega=2\omega$, we first obtain
\[
\Delta_g\rho=-\frac{2}{3\rho},
\qquad
\Delta_g(\rho\omega)=\frac{4}{3\rho}\omega.
\]
For the subsequent iterations, we have
\begin{align*}
\Delta_g(\rho^{-1})=0,\quad
\Delta_g(\rho^{-1}\omega)=2\rho^{-3}\omega,\quad
\Delta_g(\rho^{-3}\omega)=0.
\end{align*}
Direct differentiation yields
\begin{align*}
\Delta_g\psi_3
&=\left(\frac4{3\rho}\omega,-\frac{2\sqrt2}{3\rho}\right),\\
\Delta_g^2\psi_3
&=\left(\frac8{3\rho^3}\omega,0\right)\neq0,\\
\Delta_g^3\psi_3&=0.
\end{align*}
For the unit normal
$\nu=(-\sqrt{2/3}\,\omega,1/\sqrt3)$, the identity
$\Delta_g\psi_3=-3H\nu$ gives
\[
H=\frac{2\sqrt2}{3\sqrt3\,\rho}>0.
\]
Hence $\psi_3$ is nonminimal and has exact order $3$.
A radial curve has length element $\sqrt3\,\mathrm d\rho$ and reaches
the omitted vertex in finite length, so the induced metric is incomplete.
\end{example}

\begin{remark}\label{rem:ratto}
Example~\ref{ex:direct-triharmonic-cone} allows a direct comparison with
the rotational equations in \cite{Ratto2020}. In the notation of
\cite{Ratto2020}, take $p=2$, $q=0$ and
\[
x(s)=as,\qquad y(s)=bs,\qquad
a,b>0,\qquad a^2+b^2=1.
\]
The profile has unit speed and positive coordinates for $s>0$.
The corresponding rotational immersion is $\psi(s,\omega)=(as\omega,bs)$, with $\omega\in\Sph^2$ and induced
metric $g_0=\mathrm ds^2+a^2s^2g_{\Sph^2}$.
We write $\tau_3(\psi)=(\tau_{3,x}\omega,\tau_{3,y})$.

Put $\lambda=2/a^2$ and keep $g_0$ fixed. Since
$\Delta_{\Sph^2}\omega=2\omega$, we have
\[
\Delta_{g_0}(s^\alpha\omega)
=s^{\alpha-2}\bigl(\lambda-\alpha(\alpha+1)\bigr)\omega.
\]
The successive exponents $1,-1,-3$ give the factors
$\lambda-2$, $\lambda$ and $\lambda-6$, respectively.
Moreover, $\Delta_{g_0}(bs)=-2b/s$ and
$\Delta_{g_0}^2(bs)=0$. Hence $\tau_3(\psi)=-\Delta_{g_0}^3\psi$ yields
\begin{align*}
\tau_{3,x}
&=-a(\lambda-2)\lambda(\lambda-6)s^{-5}\\
&=-\frac{8(1-a^2)(1-3a^2)}{x^5},\\
\tau_{3,y}&=0.
\end{align*}
The same expression follows from the fixed-metric reduced $3$-energy
in \cite[p.~1101]{Ratto2020}.

For comparison, substitute the same profile into the printed formula
(3.4) of \cite{Ratto2020}. All terms after the first one in
$\tau_{3,x}$ vanish, since $p=2$, $q=0$ and
$x^{(j)}=y^{(j)}=0$ for $j\geq2$. Since $b^2=1-a^2$, the remaining
expression becomes
\begin{align*}
\tau_{3,x}^{\mathrm{printed}}
&=-\frac{2}{x^5}
  \bigl[-4-10a^2+2a^6+a^4(12+2b^2)\bigr]\\
&=-\frac{2}{x^5}\bigl[-4-10a^2+14a^4\bigr]
\\
&=\frac{4(1-a^2)(7a^2+2)}{x^5},
\end{align*}
where $-4-10a^2+14a^4=2(a^2-1)(7a^2+2)$.
For $a^2=1/3$, we have $x=\rho=s/\sqrt3$. The direct calculation gives
$\tau_{3,x}=0$, whereas
the printed expression gives $\tau_{3,x}^{\mathrm{printed}}=104/(9x^5)$.
Thus the printed formula (3.4) differs from the fixed-metric triharmonic
equation on the profile of Example~\ref{ex:direct-triharmonic-cone}.
\end{remark}

\begin{remark}
The restriction of \eqref{eq:seed} to $(1,2)\times\Sph^2$ has finite
higher-order energies and exact order $k$. For $n>3$, the same conclusion
holds after restricting the Euclidean factor to a bounded open set.
The restricted metrics are incomplete.
\end{remark}

\section{Circular cones and closed links}\label{sec:classification}

We first apply Lemma~\ref{lem:cone} to standard circular cones.

\begin{proof}[Proof of Theorem~\ref{thm:circular}]
Write $b=\sqrt{1-a^2}$. The metric induced by
\eqref{eq:circular-cone} is
$\mathrm ds^2+a^2s^2g_{\Sph^{n-1}}$.
The link eigenvalues of the spherical component and the last
component are $\lambda=(n-1)/a^2$ and zero, respectively. Thus
Lemma~\ref{lem:cone} gives
\begin{equation}\label{eq:general-iterates}
\Delta^p\psi_a^n
=s^{1-2p}\bigl(aP_{p,n}(\lambda)\omega,\,
bP_{p,n}(0)\bigr).
\end{equation}

Suppose first that $\psi_a^n$ is $k$-harmonic. Since $b>0$,
the last component in \eqref{eq:general-iterates} forces
$q_{j,n}=0$ for some integer $0\leq j\leq k-1$. Since $1-2j\neq0$,
the equation $q_{j,n}=0$ becomes
\[
j=\frac{n-1}{2}.
\]
Hence $n$ is odd, and the last component first vanishes after
$(n+1)/2$ iterations.

Next, the spherical component in \eqref{eq:general-iterates}
requires $\lambda=q_{j,n}$ for some $0\leq j\leq k-1$. Since $0<a<1$, we have $\lambda>n-1$.
The first identity in \eqref{eq:root-difference} gives
$j>n/2$, or $j\geq(n+1)/2$. For every $0\leq\ell<j$, the
second identity in \eqref{eq:root-difference} yields
\[
q_{j,n}-q_{\ell,n}
=2(j-\ell)\bigl(2(j+\ell)-n\bigr)>0.
\]
Hence the spherical component first vanishes at the $(j+1)$st iteration.
Since $j\geq(n+1)/2$, \eqref{eq:general-iterates} yields exact order $j+1$.

It follows that exact order $k$ requires $j=k-1$, and hence
\[
n\text{ is odd},\quad
\frac{n+1}{2}\leq k-1,\quad
\frac{n-1}{a^2}=q_{k-1,n}=(2k-3)(2k-n-1).
\]
Thus \eqref{eq:aperture} holds.

Conversely, assume \eqref{eq:aperture}. Since $3\leq n\leq2k-3$,
\[
q_{k-1,n}-(n-1)=2(k-1)(2k-n-2)>0.
\]
Thus $a^2=(n-1)/q_{k-1,n}$ lies in $(0,1)$.
As above, the last component vanishes after $(n+1)/2\leq k-1$
iterations, whereas the spherical component first vanishes after $k$
iterations by the strict root inequalities above. Thus
\eqref{eq:general-iterates} gives exact order $k$.
\end{proof}

The principal curvatures of \eqref{eq:circular-cone} are zero once
and $b/(as)$ with multiplicity $n-1$. Under the assumptions of
Theorem~\ref{thm:circular}, substituting \eqref{eq:aperture}, we obtain
\[
|A|^2=\frac{2(k-1)(2k-n-2)}{s^2},
\qquad
H^2=\frac{2(k-1)(n-1)(2k-n-2)}{n^2s^2}.
\]
In particular, the mean curvature and the squared norm of the second
fundamental form are nonconstant.

We next consider cones over arbitrary closed spherical links.

\begin{proof}[Proof of Theorem~\ref{thm:closed-link}]
Put $d=n-1$. Since $\psi_L$ is isometric and $|\psi_L|=1$, the induced
cone metric is $g=\mathrm ds^2+s^2h$.
By Lemma~\ref{lem:cone}, the $k$-harmonic equation becomes
\begin{equation}\label{eq:link-polynomial}
P_{k,n}(\Delta_h)\psi_L=0.
\end{equation}
Since $L$ is closed, $\Delta_h$ is nonnegative and self-adjoint
with discrete spectrum. If $\psi_{L,\lambda}$ denotes the orthogonal
projection of $\psi_L$ onto the $\lambda$-eigenspace, then
\eqref{eq:link-polynomial} gives
$P_{k,n}(\lambda)\psi_{L,\lambda}=0$. Consequently,
\[
\psi_L=\sum_\lambda\psi_{L,\lambda},
\qquad
\Delta_h\psi_{L,\lambda}=\lambda\psi_{L,\lambda},
\]
where the finite sum runs over the distinct nonnegative roots of
$P_{k,n}$ that are eigenvalues of $\Delta_h$. This conclusion also
holds when $P_{k,n}$ has repeated roots.

We now examine the positive roots of $P_{k,n}$. We have $q_{0,n}=d$.
For $j\geq1$, the conditions $q_{j,n}>0$
and $1-2j<0$ yield $d-2j<0$. It follows that $2j>d$, and hence $2j\geq n$.
Using \eqref{eq:root-difference}, we obtain
\begin{equation}\label{eq:positive-root-bound}
q_{j,n}-d=2j(2j-n)\geq0.
\end{equation}

Since $|\psi_L|^2=1$ and $|\mathrm d\psi_L|^2=d$, orthogonality and
integration by parts give
\begin{equation}\label{eq:link-norms}
\sum_\lambda\|\psi_{L,\lambda}\|_{L^2}^2=\vol(L),
\qquad
\sum_\lambda\lambda\|\psi_{L,\lambda}\|_{L^2}^2=d\,\vol(L).
\end{equation}
By \eqref{eq:positive-root-bound}, every positive root is at least $d$.
Subtracting $d$ times the first identity from the second therefore yields
\begin{equation}\label{eq:spectral-identity}
\sum_{\lambda>d}(\lambda-d)\|\psi_{L,\lambda}\|_{L^2}^2
=d\|\psi_{L,0}\|_{L^2}^2.
\end{equation}
Here $\psi_{L,0}=0$ if zero does not occur in the spectral decomposition.

We claim that $\psi_{L,0}\neq0$. Indeed, if $\psi_{L,0}=0$, then
\eqref{eq:spectral-identity} forces every
component with $\lambda>d$ to vanish. Thus $\Delta_h\psi_L=d\psi_L$,
and \eqref{eq:homogeneous} gives $\Delta_g\psi=0$, which contradicts
nonminimality. Hence $\psi_{L,0}\neq0$, and
\eqref{eq:spectral-identity} also forces a nonzero component with
$\lambda>d$.

The occurrence of the zero eigenvalue implies that $q_{j,n}=0$
for some integer $j$. As in the proof of
Theorem~\ref{thm:circular}, $n$ must be odd and $j=(n-1)/2$.
The first root strictly greater than $d$ has index $(n+1)/2$.
Hence
\[
\frac{n+1}{2}\leq k-1,
\qquad\text{and therefore}\qquad
3\leq n\leq2k-3.
\]
Since $L$ is connected and closed, every harmonic function on $L$ is
constant. Hence
\[
\psi_{L,0}=\frac1{\vol(L)}\int_L\psi_L\dd v_h\neq0.
\]
Finally, Theorem~\ref{thm:circular} gives a cone of exact order $k$
with a closed link for every permitted pair $(k,n)$.
\end{proof}

\begin{proof}[Proof of Corollary~\ref{cor:centroid}]
Put $d=n-1$. By Theorem~\ref{thm:closed-link}, $n$ is odd,
$3\leq n\leq2k-3$, and $c\neq0$. We have
$\Lambda=q_{k-1,n}$. For every $0\leq j\leq k-2$,
\eqref{eq:root-difference} yields
\[
\Lambda-q_{j,n}
=2(k-1-j)\bigl(2(k-1+j)-n\bigr)>0,
\]
since $2(k-1)-n\geq1$. Hence $\Lambda$ is the largest root of
$P_{k,n}$ and $\Lambda>d$.

Consider the spectral decomposition in the proof of
Theorem~\ref{thm:closed-link}. Since $L$ is connected and closed, the
zero spectral component is the constant vector $c$, and therefore
\[
\psi_L=c+\sum_{\lambda>0}\psi_{L,\lambda},
\qquad
\|\psi_{L,0}\|_{L^2}^2=|c|^2\vol(L).
\]
Every positive eigenvalue occurring in this sum is at most $\Lambda$.
Combining the spectral decomposition with \eqref{eq:link-norms},
we obtain
\[
\begin{aligned}
d\,\vol(L)
&=\sum_{\lambda>0}\lambda\|\psi_{L,\lambda}\|_{L^2}^2\\
&\leq\Lambda\sum_{\lambda>0}\|\psi_{L,\lambda}\|_{L^2}^2\\
&=\Lambda(1-|c|^2)\vol(L).
\end{aligned}
\]
Since $c\neq0$, this yields \eqref{eq:centroid-bound}.

Suppose that equality holds. Then
\[
\sum_{0<\lambda<\Lambda}
(\Lambda-\lambda)\|\psi_{L,\lambda}\|_{L^2}^2=0.
\]
Every coefficient is positive, so every positive spectral component
with $\lambda<\Lambda$ vanishes. Consequently,
\[
\psi_L=c+\psi_{L,\Lambda},
\qquad
\Delta_h\psi_L=\Lambda(\psi_L-c).
\]
The identities $|\psi_L|^2=1$ and $|\mathrm d\psi_L|^2=d$ now give
\[
\begin{aligned}
0=\frac12\Delta_h|\psi_L|^2
&=\langle\Delta_h\psi_L,\psi_L\rangle
  -|\mathrm d\psi_L|^2\\
&=\Lambda\bigl(1-\langle c,\psi_L\rangle\bigr)-d.
\end{aligned}
\]
Thus $\langle c,\psi_L\rangle=1-d/\Lambda$ pointwise. Integrating
over $L$ and using the definition of $c$, we obtain
\[
\langle c,\psi_L\rangle=|c|^2=1-\frac d\Lambda.
\]
Put $a=\sqrt{d/\Lambda}$. Since $\Lambda>d$, we have $0<a<1$.
It follows that
\[
\langle c,\psi_L-c\rangle=0,
\qquad
|\psi_L-c|^2=a^2.
\]
Hence $\widehat\psi=a^{-1}(\psi_L-c)$ takes values in the unit sphere
of $c^\perp$, which is identified with $\Sph^{N-2}(1)$. Since $c$ is
constant and $\psi_L$ is isometric,
\[
\widehat\psi^{\,*}g_{\Sph^{N-2}}=a^{-2}h.
\]
Under this constant rescaling of the metric,
\[
\Delta_{a^{-2}h}\widehat\psi
=a^2\Delta_h\widehat\psi
=a^2\Lambda\widehat\psi
=d\widehat\psi.
\]
Takahashi's characterization \cite{Takahashi1966} therefore yields
minimality of $\widehat\psi$. An orthogonal transformation sending
$c$ to $\sqrt{1-a^2}\,e_N$ gives \eqref{eq:centroid-equality}.

Conversely, suppose that \eqref{eq:centroid-equality} holds with the
stated minimal isometric immersion $\widehat\psi$. Write
$b=\sqrt{1-a^2}$. Takahashi's characterization and metric rescaling
give
\[
\Delta_h\widehat\psi
=a^{-2}\Delta_{a^{-2}h}\widehat\psi
=\frac d{a^2}\widehat\psi
=\Lambda\widehat\psi.
\]
Since $L$ is closed, integration of this identity gives
\[
\Lambda\int_L\widehat\psi\dd v_h
=\int_L\Delta_h\widehat\psi\dd v_h=0.
\]
Thus $c=(0,b)$ and $|c|^2=b^2=1-d/\Lambda$, which proves equality
in \eqref{eq:centroid-bound}.

Finally, $\Lambda=q_{k-1,n}$ and
$q_{(n-1)/2,n}=0$, where $(n-1)/2\leq k-2$. Thus
\[
P_{k,n}(\Lambda)=P_{k,n}(0)=P_{k-1,n}(0)=0,
\qquad
P_{k-1,n}(\Lambda)
=\prod_{j=0}^{k-2}(\Lambda-q_{j,n})>0.
\]
Applying Lemma~\ref{lem:cone} to
$\psi(s,x)=s(a\widehat\psi(x),b)$, we obtain
\[
\Delta_g^k\psi=0,
\qquad
\Delta_g^{k-1}\psi
=s^{3-2k}\bigl(aP_{k-1,n}(\Lambda)\widehat\psi,0\bigr)\neq0.
\]
The last inequality holds at every point because $|\widehat\psi|=1$.
Hence $\psi$ has exact order $k$. Taking
$L=\Sph^{n-1}$, $h=a^2g_{\Sph^{n-1}}$ and
$\widehat\psi=\operatorname{id}_{\Sph^{n-1}}$ gives the circular
cones in $\R^{n+1}$ from Theorem~\ref{thm:circular}, so the upper
bound is sharp.
\end{proof}

\begin{remark}
For $n>3$, the smooth link $M_{k,n}\cap\Sph^n(1)$ is diffeomorphic to
$\Sph^2\times B^{n-3}$ through the parametrization
\[
(\omega,z)\longmapsto
\bigl(a_k\sqrt{1-|z|^2}\,\omega,\ z,\,
b_k\sqrt{1-|z|^2}\bigr),
\]
where $B^{n-3}$ is the open unit ball in $\R^{n-3}$.
Thus this link is noncompact, and the even-dimensional examples
in Theorem~\ref{thm:main} do not fall under
Theorem~\ref{thm:closed-link}.
\end{remark}

\begin{remark}
The finite spectral decomposition above concerns the link map
$\psi_L$ and the link Laplacian $\Delta_h$. It does not give a finite
spectral decomposition of the cone immersion with respect to
$\Delta_g$. Indeed, suppose that $\psi$ is a finite sum of
$\Delta_g$-eigenmaps and satisfies $\Delta_g^k\psi=0$.
Since distinct eigenspaces are linearly independent, every component
with nonzero eigenvalue vanishes. Hence $\Delta_g\psi=0$.
The nonminimal cones above are therefore not of finite type.
This is a special case of \cite[Proposition~3.1]{Chen2007}.
For related results on finite-type cones and ruled submanifolds,
see \cite{Garay1988,FGL1991}.
\end{remark}

\section{Examples in higher codimension}\label{sec:examples}

We now extend the construction using minimal immersions into spheres.
Let $\psi_L:(L^d,h)\to\Sph^{N-1}(1)$ be a minimal isometric immersion.
Takahashi's characterization \cite{Takahashi1966} gives
$\Delta_h\psi_L=d\psi_L$.

\begin{proposition}\label{prop:blocks}
Let $k\geq3$ and let $n$ be odd with $3\leq n\leq2k-3$.
Put $d=n-1$. For $i=1,\ldots,J$, let
\[
\psi_i:(L_i^{d_i},h_i)\longrightarrow\Sph^{N_i-1}(1)
\]
be minimal isometric immersions, where $d_i\geq1$ and
$\sum_i d_i=d$. Choose integers
\[
j_i\in\{0\}\cup\left\{\frac{n+1}{2},\ldots,k-1\right\},
\qquad
\max_i j_i=k-1,
\]
and define
\begin{equation}\label{eq:block-parameters}
\lambda_i=q_{j_i,n},
\qquad
a_i=\sqrt{\frac{d_i}{\lambda_i}},
\qquad
b=\sqrt{1-\sum_{i=1}^J a_i^2}.
\end{equation}
Then $b>0$, and the map
\begin{equation}\label{eq:block-cone}
\psi(s,x_1,\ldots,x_J)
=s\bigl(a_1\psi_1(x_1),\ldots,a_J\psi_J(x_J),b\bigr)
\end{equation}
from $(0,\infty)\times\prod_iL_i$ into $\R^{1+\sum_iN_i}$
is a nonminimal immersion of exact order $k$. Moreover,
\begin{equation}\label{eq:block-mean}
|\Hm|^2
=\frac{\sum_{i=1}^J d_i\lambda_i-d^2}{n^2s^2}>0.
\end{equation}
\end{proposition}

\begin{proof}
First, the choice of the indices and \eqref{eq:positive-root-bound}
give $\lambda_i\geq d$ for every $i$, with strict inequality for at
least one index. Hence
\begin{equation*}
\sum_i a_i^2=\sum_i\frac{d_i}{\lambda_i}
<\frac{\sum_i d_i}{d}=1.
\end{equation*}
It follows that $b>0$. Moreover, the radial derivative of
\eqref{eq:block-cone} has length one and is orthogonal to every
derivative tangent to the factors. Since each $\psi_i$ is isometric,
the induced metric has the form
$$
g=\mathrm ds^2+s^2h,
\qquad
h=\sum_i a_i^2h_i.
$$
In particular, $\psi$ is an immersion.

Next, the equality $\Delta_{h_i}\psi_i=d_i\psi_i$ gives
$\Delta_h\psi_i=\lambda_i\psi_i$. The last component is constant on
the link. Applying Lemma~\ref{lem:cone}, we obtain
\begin{equation}\label{eq:block-iterates}
\Delta_g^p\psi
=s^{1-2p}
\bigl(a_1P_{p,n}(\lambda_1)\psi_1,\ldots,
a_JP_{p,n}(\lambda_J)\psi_J,bP_{p,n}(0)\bigr).
\end{equation}
If $j_i=0$, then $\lambda_i=q_{0,n}=d$, and the $i$th component
vanishes at the first iteration. If $j_i\geq(n+1)/2$, then
\eqref{eq:root-difference} gives
\[
\lambda_i-q_{\ell,n}
=2(j_i-\ell)\bigl(2(j_i+\ell)-n\bigr)>0
\qquad(0\leq\ell<j_i).
\]
Hence the $i$th component first vanishes after $j_i+1$ iterations.
The last component first vanishes after $(n+1)/2\leq k-1$
iterations, since $q_{(n-1)/2,n}=0$.
For an index $i$ with $j_i=k-1$, the strict inequalities above give
$P_{k-1,n}(\lambda_i)>0$. Since $|\psi_i|=1$,
\eqref{eq:block-iterates} shows that the $i$th component of
$\Delta_g^{k-1}\psi$ is nowhere zero. Hence
$\Delta_g^k\psi=0$ and $\Delta_g^{k-1}\psi\neq0$ everywhere.

Finally, taking $p=1$ in \eqref{eq:block-iterates}, we obtain
\[
\Delta_g\psi
=s^{-1}\bigl(a_1(\lambda_1-d)\psi_1,\ldots,
a_J(\lambda_J-d)\psi_J,-db\bigr).
\]
By \eqref{eq:block-parameters}, $a_i^2\lambda_i=d_i$ and
$\sum_i a_i^2+b^2=1$. Since $\Delta_g\psi=-n\Hm$, we obtain
\[
\begin{aligned}
n^2s^2|\Hm|^2
&=\sum_i a_i^2(\lambda_i-d)^2+d^2b^2\\
&=\sum_i a_i^2\lambda_i^2-2d\sum_i a_i^2\lambda_i
  +d^2\left(\sum_i a_i^2+b^2\right)\\
&=\sum_i d_i\lambda_i-d^2.
\end{aligned}
\]
Since at least one $\lambda_i>d$, we obtain \eqref{eq:block-mean}.
\end{proof}

Taking $J=1$ in Proposition~\ref{prop:blocks} gives
\[
\psi(s,x)=(as\widehat\psi(x),bs),
\qquad
a^2=\frac{n-1}{(2k-3)(2k-n-1)},
\qquad b=\sqrt{1-a^2},
\]
where $\widehat\psi:(L^{n-1},h)\to\Sph^{N-1}(1)$ is any minimal
isometric immersion. If $L$ is connected and closed, integration of
$\Delta_h\widehat\psi=(n-1)\widehat\psi$ gives
$\int_L\widehat\psi\dd v_h=0$. The cone link has metric $a^2h$ and
average position $(0,b)$, and hence attains equality in
\eqref{eq:centroid-bound}. Conversely, Corollary~\ref{cor:centroid}
shows that every equality case arises in this way, up to an ambient
orthogonal transformation. For the unit normal field
$\nu=(-b\widehat\psi,a)$, direct calculation yields
\[
\Hm=\frac{(n-1)b}{nas}\nu.
\]
If the constant component $b$ of the link is omitted and
$\sum_i a_i^2=1$, then any $k$-harmonic cone constructed from minimal
spherical factors must be minimal. Indeed, \eqref{eq:positive-root-bound}
yields $\lambda_i\geq d$, and hence
\[
1=\sum_i\frac{d_i}{\lambda_i}
\leq\sum_i\frac{d_i}{d}=1.
\]
Equality forces every $\lambda_i=d$, and therefore $\Delta_g\psi=0$.

To give explicit examples, put $e(t)=(\cos t,\sin t)$ and take
$s>0$ and $\theta_1,\theta_2\in\R/(2\pi\mathbb Z)$.

\begin{figure}[htbp]
\centering
\includegraphics[width=0.72\textwidth]{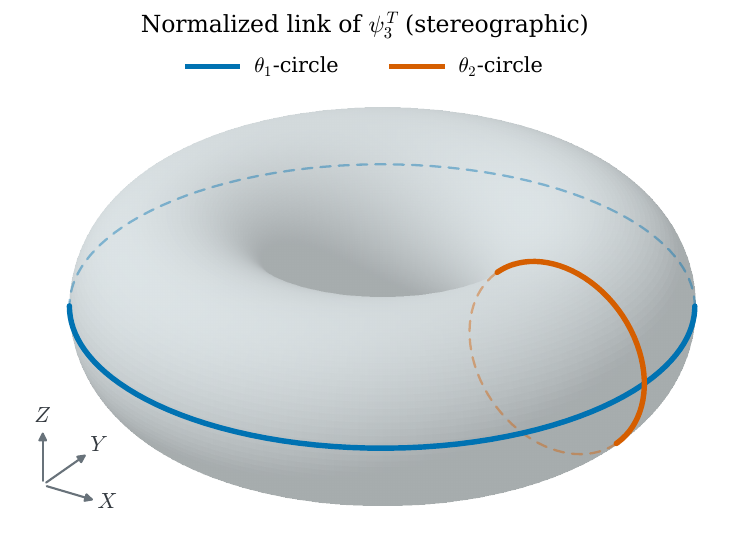}
\caption{Stereographic image of the normalized link of $\psi_3^{\mathrm T}$.
The projection is taken from $(0,0,0,1)$.
The blue and orange curves are the $\theta_1$- and
$\theta_2$-coordinate circles, respectively.}
\label{fig:torus-link}
\end{figure}

\begin{example}
For every $k\geq3$, set $\Lambda_k=2(k-2)(2k-3)$ and define
\begin{equation}\label{eq:torus-family}
\psi_k^{\mathrm T}(s,\theta_1,\theta_2)
=
\left(
\frac{s}{\sqrt2}e(\theta_1),\
\frac{s}{\sqrt{\Lambda_k}}e(\theta_2),\
s\sqrt{\frac12-\frac1{\Lambda_k}}
\right)\in\R^5.
\end{equation}
The induced metric is given by
\[
g=\mathrm ds^2+\frac{s^2}{2}\,\mathrm d\theta_1^2
+\frac{s^2}{\Lambda_k}\,\mathrm d\theta_2^2.
\]
The three components have link eigenvalues $2$, $\Lambda_k$,
and zero. Since $\Lambda_k>2$, two distinct positive spectral values
occur, so the upper bound in \eqref{eq:centroid-bound} is strict.
Proposition~\ref{prop:blocks} gives exact order $k$ and
\[
|\Hm|^2=\frac{\Lambda_k-2}{9s^2}>0.
\]
Since the link metric is flat, the cone has scalar curvature $-2/s^2$.
Combining this with the Gauss equation $\Scal=9|\Hm|^2-|\II|^2$, we obtain
\[
|\II|^2=\frac{\Lambda_k}{s^2}.
\]
For $k=3$, \eqref{eq:torus-family} becomes
\[
\psi_3^{\mathrm T}=
\left(
\frac{s}{\sqrt2}e(\theta_1),\
\frac{s}{\sqrt6}e(\theta_2),\
\frac{s}{\sqrt3}
\right).
\]
Since $(q_{0,3},q_{1,3},q_{2,3})=(2,0,6)$,
Lemma~\ref{lem:cone} gives
\[
\Delta^2\psi_3^{\mathrm T}
=\left(0,\frac{24}{\sqrt6\,s^3}e(\theta_2),0\right)\neq0,
\qquad
\Delta^3\psi_3^{\mathrm T}=0.
\]
We also obtain $|\Hm|^2=4/(9s^2)$ and $|\II|^2=6/s^2$.
Figure~\ref{fig:torus-link} shows the stereographic image of
$(\sqrt3\,e(\theta_1)/2,e(\theta_2)/2)\in\Sph^3$, obtained by normalizing
the first four coordinates of the link $s=1$.
\end{example}

\begin{example}
Consider
\[
\psi_4^{\mathrm U}(s,\theta_1,\theta_2)
=
\left(
\frac{s}{\sqrt6}e(\theta_1),\
\frac{s}{\sqrt{20}}e(\theta_2),\
\sqrt{\frac{47}{60}}\,s
\right)\in\R^5.
\]
The induced metric is given by
\[
g=\mathrm ds^2+\frac{s^2}{6}\,\mathrm d\theta_1^2
+\frac{s^2}{20}\,\mathrm d\theta_2^2.
\]
We have $(q_{0,3},q_{1,3},q_{2,3},q_{3,3})=(2,0,6,20)$.
Here both positive link eigenvalues, $6$ and $20$, exceed $n-1=2$.
The first spherical component vanishes after three iterations,
the second after four, and the last component after two.
Applying Lemma~\ref{lem:cone}, we obtain
\[
\Delta^3\psi_4^{\mathrm U}
=\left(0,\frac{5040}{\sqrt{20}\,s^5}e(\theta_2),0\right)\neq0,
\qquad
\Delta^4\psi_4^{\mathrm U}=0.
\]
Moreover, equation \eqref{eq:block-mean} and the Gauss equation yield
\[
|\Hm|^2=\frac{22}{9s^2},
\qquad
|\II|^2=\frac{24}{s^2}.
\]
In particular, the length of the mean curvature vector is nonconstant.
\end{example}
\begin{remark}
	All the examples constructed in this paper are nonminimal and incomplete.
	Indeed, the omitted cone vertices can be approached along radial curves
	of finite length, and this incompleteness persists under taking Euclidean
	products. Thus, the counterexamples above leave open the following
	complete version of the $k$-harmonic conjecture.
	
	\smallskip
	\noindent
\textbf{ $\boldsymbol{k}$-harmonic conjecture of complete version.}
	\emph{For each integer $k\geq2$, every complete $k$-harmonic submanifold
		of Euclidean space is minimal.}
	
	\smallskip
	For $k=2$, this is Chen's biharmonic conjecture restricted to complete
	submanifolds. The minimality results under additional integrability
	assumptions in \cite{MNU2015,NU2018} provide partial affirmative answers.
	It remains to determine whether completeness alone is sufficient
	to force minimality.
\end{remark}

\section*{Acknowledgement}
The work was partially supported by NSFC (No.~12571051). The first author wish to thank Professor Cezar Oniciuc for his useful comments, which helped to improve the original version of this paper.

\section*{Declaration of generative AI}

The authors used AI to assist with language editing and technical checks. The authors take full responsibility for the content of this manuscript.


\begin{thebibliography}{99}

\bibitem{AkutagawaMaeta2013}
K.~Akutagawa, S.~Maeta,
Biharmonic properly immersed submanifolds in Euclidean spaces,
Geom. Dedicata 164 (2013), 351--355.

\bibitem{BMO2008}
A.~Balmu\c{s}, S.~Montaldo, C.~Oniciuc,
Classification results for biharmonic submanifolds in spheres,
Israel J. Math. 168 (2008), 201--220.

\bibitem{BMO2010}
A.~Balmu\c{s}, S.~Montaldo, C.~Oniciuc,
Biharmonic hypersurfaces in 4-dimensional space forms,
Math. Nachr. 283 (12) (2010), 1696--1705.

\bibitem{Bernstein1915}
S.~Bernstein,
Sur un th\'eor\`eme de g\'eom\'etrie et son application aux
\'equations aux d\'eriv\'ees partielles du type elliptique,
Comm. Soc. Math. Kharkow (2) 15 (1) (1915),
38--45.

\bibitem{BDG1969}
E.~Bombieri, E.~De Giorgi, E.~Giusti,
Minimal cones and the Bernstein problem,
Invent. Math. 7 (1969), 243--268.

\bibitem{Branding2023}
V.~Branding,
On the normal stability of triharmonic hypersurfaces in space forms,
J. Geom. Anal. 33 (11) (2023), 355.

\bibitem{Branding2025}
V.~Branding,
On the normal stability of the 4-harmonic and the ES-4-harmonic hypersphere,
J. Differential Equations 424 (2025), 586--636.

\bibitem{BMOR2020}
V.~Branding, S.~Montaldo, C.~Oniciuc, A.~Ratto,
Higher order energy functionals,
Adv. Math. 370 (2020), 107236.

\bibitem{Chen1991}
B.-Y.~Chen,
Some open problems and conjectures on submanifolds of finite type,
Soochow J. Math. 17 (2) (1991), 169--188.

\bibitem{Chen2007}
B.-Y.~Chen,
Tension field, iterated Laplacian, type number and Gauss maps,
Houston J. Math. 33 (2) (2007), 461--481.

\bibitem{ChenGuan2023}
H.~Chen, Z.~Guan,
Triharmonic CMC hypersurfaces in $\R^5(c)$,
Manuscripta Math. 172 (1--2) (2023), 209--220.

\bibitem{FGL1991}
A.~Ferr\'andez, O.~J.~Garay, P.~Lucas,
Finite type ruled manifolds shaped on spherical submanifolds,
Arch. Math. (Basel) 57 (1) (1991), 97--104.

\bibitem{FuHongTianJEMS}
Y.~Fu, M.-C.~Hong, G.~Tian,
The biharmonic hypersurface flow and the Willmore flow in higher dimensions,
J. Eur. Math. Soc. (2026), online first, doi:10.4171/JEMS/1770.

\bibitem{FHY2026}
Y.~Fu, M.-C.~Hong, D.~Yang,
Biharmonic hypersurfaces in space forms, arXiv:2610.04165v2, 2026.

\bibitem{FHZ2021}
Y.~Fu, M.-C.~Hong, X.~Zhan,
On Chen's biharmonic conjecture for hypersurfaces in $\R^5$,
Adv. Math. 383 (2021), 107697.

\bibitem{FHZ2023}
Y.~Fu, M.-C.~Hong, X.~Zhan,
Biharmonic conjectures on hypersurfaces in a space form,
Trans. Amer. Math. Soc. 376 (12) (2023),
8411--8445.

\bibitem{FuYang2023}
Y.~Fu, D.~Yang,
On triharmonic hypersurfaces in space forms,
J. Geom. Anal. 33 (8) (2023), 263.

\bibitem{Garay1988}
O.~J.~Garay,
Finite type cones shaped on spherical submanifolds,
Proc. Amer. Math. Soc. 104 (1988), 868--870.

\bibitem{HasanisVlachos1995}
Th.~Hasanis, Th.~Vlachos,
Hypersurfaces in $E^4$ with harmonic mean curvature vector field,
Math. Nachr. 172 (1995), 145--169.

\bibitem{Jiang1986}
G.~Y.~Jiang,
2-harmonic maps and their first and second variational formulas,
Chinese Ann. Math. Ser. A 7 (4) (1986),
389--402 (in Chinese).

\bibitem{Jiang1987}
G.~Y.~Jiang,
Some nonexistence theorems on 2-harmonic and isometric immersions
in Euclidean space,
Chinese Ann. Math. Ser. A 8 (3) (1987),
377--383 (in Chinese).

\bibitem{Maeta2012}
S.~Maeta,
$k$-harmonic maps into a Riemannian manifold with constant sectional curvature,
Proc. Amer. Math. Soc. 140 (5) (2012), 1835--1847.
Preprint: arXiv:1010.0920v1.

\bibitem{Maeta2015Construction}
S.~Maeta,
Construction of triharmonic maps,
Houston J. Math. 41 (2) (2015), 433--444.

\bibitem{MNU2015}
S.~Maeta, N.~Nakauchi, H.~Urakawa,
Triharmonic isometric immersions into a manifold of non-positively constant
curvature,
Monatsh. Math. 177 (4) (2015), 551--567.

\bibitem{MOR2022}
S.~Montaldo, C.~Oniciuc, A.~Ratto,
Polyharmonic hypersurfaces into space forms,
Israel J. Math. 249 (1) (2022), 343--374.

\bibitem{MOR2024}
S.~Montaldo, C.~Oniciuc, A.~Ratto,
Polyharmonic surfaces in 3-dimensional homogeneous spaces,
Manuscripta Math. 174 (2024), 31--58.

\bibitem{MR2018}
S.~Montaldo, A.~Ratto,
New examples of $r$-harmonic immersions into the sphere,
J. Math. Anal. Appl. 458 (1) (2018), 849--859.

\bibitem{NU2018}
N.~Nakauchi, H.~Urakawa,
Polyharmonic maps into the Euclidean space,
Note Mat. 38 (1) (2018), 89--100.


\bibitem{Ratto2020}
A.~Ratto,
Higher order energy functionals and the Chen--Maeta conjecture,
AIMS Math. 5 (2) (2020), 1089--1104.

\bibitem{Simons1968}
J.~Simons,
Minimal varieties in Riemannian manifolds,
Ann. of Math. (2) 88 (1) (1968), 62--105.

\bibitem{Takahashi1966}
T.~Takahashi,
Minimal immersions of Riemannian manifolds,
J. Math. Soc. Japan 18 (4) (1966), 380--385.

\end{thebibliography}
\end{document}